\documentclass[fleqn,10pt]{wlscirep}%10pt
\usepackage[utf8]{inputenc}
\usepackage[T1]{fontenc}
\usepackage{amsmath}
\usepackage{graphicx} 
\usepackage{paralist, tabularx}
\usepackage[outercaption]{sidecap}
\title{\LARGE{Sounding the metabolic orchestra:\\ A delay dynamical systems perspective on the glucose-insulin regulatory response to on-off glucose infusion}}
\author[1,2,*]{Stefan Ruschel}
\author[1,+]{Benoit Huard}

\affil[1]{Northumbria University, Department of Mathematics, Physics and Electrical Engineering, Newcastle upon Tyne, NE1 8ST, United Kingdom}
\affil[2]{The University of Auckland, Department of Mathematics, Auckland, 1142, New Zealand}

\affil[*]{stefan.ruschel@northumbria.ac.uk}
\affil[+]{benoit.huard@northumbria.ac.uk}

\begin{abstract}
We investigate the consequences of periodic, on-off glucose infusion on the glucose-insulin regulatory system on the basis of a system-level mathematical model with two explicit time delays. 
Studying the effects of such infusion protocols is mathematically challenging yet a promising direction for probing the system response to infusion. 
We pay special attention to the interplay of the infusion with intermediate-time-scale, ultradian oscillations that arise as a results of the physiological response of glucose uptake and back-release into the bloodstream.
By using numerical solvers and numerical continuation software, we investigate the response of the model to different infusion patterns, and explore how these patterns affect the overall levels of glucose and insulin, and can lead to entrainment. 
By doing so, we provide a road-map of system responses that can potentially help identify new test strategies for detecting abnormal responses to glucose uptake.
\end{abstract}

\begin{document}

\flushbottom
\maketitle
\thispagestyle{empty}

\section*{Introduction}

%Introduce the concept of cyclic rhythms and their role in regulating biological and physiological systems. Example are blood glucose levels
Cyclic rhythms are widely recognized for their significant role in regulating the function of biological and physiological systems
\cite{goldbeter2022multi, keener2009mathematical}. 
Endogenic oscillations are typically encountered in healthy individuals, while a progressive lack of control of these rhythms is often associated with system stress (e.g. sleep deprivation \cite{sweeney2017skeletal, sweeney2021impairments}), and disease evolution in humans \cite{spiga2011hpa}. 

A prominent example of such endocrine oscillations in the human body is the self-regulation of blood glucose levels \cite{grant2018evidence}. 
When blood glucose levels increase, insulin is released from the pancreas. 
Insulin then causes blood glucose levels to decrease by stimulating body cells to absorb glucose from the blood. 
Conversely, when blood glucose levels fall, pancreatic $\beta$-cells release glucagon stimulating hepatic glycogenolysis and neoglucogenesis. 
The level of blood glucose is then controlled by the rates of insulin secretion (activation by glucose) and hepatic glucose production (inhibition by insulin).
Within the glucose-insulin regulatory system, both rapid oscillations of insulin (period $\sim$ 6-15 minutes), and ultradian oscillation of glucose and insulin (of similar period  $\sim$ 80-180 minutes \cite{scheen1996relationships}) have been observed during fasting, meal ingestion, continuous enteral and intravenous nutrition \cite{o1993lack}.

The most important pathway to understand the underlying mechanisms of these glucose-insulin oscillations is measuring the response to glucose infusions.
A large quantity of metrics and mathematical models have been devised for that purpose. 
While the HBA1c metric remains an essential tool for the diagnostic, prevention and control for T2 diabetes \cite{ADA-A1C-2022}, clinical tests involving patterns of glucose intake combined with mathematical models provide a mechanism for evaluating the efficacy of internal regulation \cite{ajmera2013impact,makroglou2006mathematical,palumbo2013mathematical, HuardKirkham2022mathematical}.
The minimal model devised by Bergman and Cobelli \cite{bergman1979quantitative,bergman2021origins} provides an effective method for estimating insulin sensitivity from an intravenous or oral glucose tolerance test, although it can lead to underestimation in individuals with a large acute insulin response \cite{ha2021minmod}.
With the wider availability of continuous glucose monitors and automated insulin pumps, the ability to detect diabetic deficiencies relies on the capacity of models to reproduce more complex and realistic dynamics under various routine life conditions such as, for example, sleep deprivation \cite{sweeney2017skeletal}.

The main goal of this article is to identify the types of behaviors in a suitable mathematical model that can be expected as a response to periodic glucose uptake, specifically periodic on-off glucose infusion, which can be readily implemented in practice. 
We focus on the capacity of the system to fall into lockstep with the frequency of the glucose stimulus (so-called entrainment) which has been observed in numerous contexts at the ultradian and circadian levels in endocrinology \cite{walker2010origin, zavala2019mathematical, HuardKirkham2022mathematical}, but especially in models of glucose-insulin oscillations with periodic infusion \cite{sturis1995phase}. 

% Introduce the use of time delays in modeling oscillatory behavior in complex biological networks 
Many modeling efforts have been made to replicate the nonlinear response of the glycolytic system; in particular, the mathematical modeling of the delayed response of individual parts of the system by explicit time delays has proven an effective means to explain the onset of self-sustained, ultradian oscillations in the glycemic system \cite{li2006modeling, chen2010modeling, cohen2021novel}. 
A common approach to modeling oscillatory behavior of complex biology is to consider time delays \cite{glass2021nonlinear}.
In particular, models of endocrine regulation often incorporate explicit delays to account for the time required for the synthesis, release, and action of hormones or metabolites \cite{walker2010origin}. 
Various models of intrapancreatic rhythmic activity have been proposed recently, see Ref.~\citeonline{HuardKirkham2022mathematical} for a review. 
For example, it was shown that glucose oscillations can enhance the insulin secretory response at the $\beta$-cell level when tweaked at a suitable amplitude and frequency \cite{mckenna2016glucose}.
Negative delayed feedback has also been shown to provide a suitable explanatory mechanism for the coordinated pancreatic islet activity
%through modeling and microfluidics 
\cite{bruce_coordination_2022}.

\begin{figure}[!]
    \begin{flushleft}
    \includegraphics[width=\linewidth]{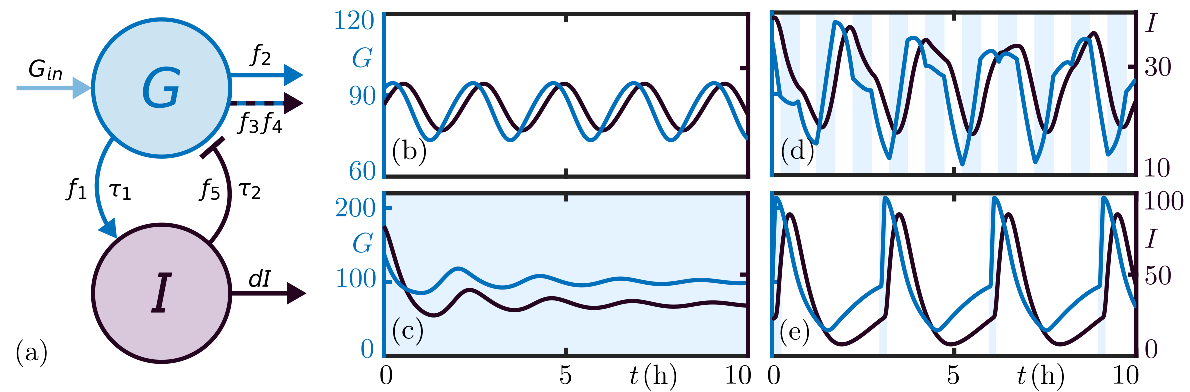}
    \end{flushleft}
    \caption{
Panel (a): Diagrammatic overview of the glucose-insulin regulatory delayed-feedback model (\ref{eq:G})--(\ref{eq:I}); see methods section for details. 
Panels (b)--(e): Characteristic time series of system (\ref{eq:G})--(\ref{eq:I}) (with positive initial condition) for different patterns of glucose infusion with intervals of fasting indicated by a white background and intervals of glucose infusion with constant rate indicated by a light blue background. 
Units are [$G$] mg dl$^{-1}$, [$I$]  mg dl$^{-1}$, and [$t$] h. Infusion rates when not fasting are $G_{\text{in}}=1.35$ mg dl$^{-1}$min$^{-1}$ in panels (c)--(d) and $G_{\text{in}}=24.3$ mg dl$^{-1}$min$^{-1}$ in panel (e); period of infusion is $T_{\text{in}}= 1$ h in panel (d) and $T_{\text{in}}= 3$ h in panel (e); time of infusion is $t_{\text{in}}= 30$ min in panel (d) and $t_{\text{in}}= 5$ min in panel (e). }
    \label{fig:schema}
\end{figure}

% Our model
%\subsection*{The glucose-insulin regulatory delayed feedback model}
In this article, we investigate a two-component, system-level mathematical model (see Eqs.~\eqref{eq:G}--\eqref{eq:I} in the methods section) for blood glucose level $G(t)$ and insulin levels $I(t)$ with two explicit time delays $\tau_I$ and $\tau_G$
corresponding to pancreatic insulin and hepatic glucose production pathways, see the methods section for details on the model. 
%%%%%%%%%%%%%%%%%%%%%%%%%%%
The model incorporates the following physiological processes and factors that influence glucose and insulin dynamics, see Fig.~\ref{fig:schema}(a) for a schematic overview.
\begin{compactitem}
   \item Glucose uptake: $G_\text{in}$ represents glucose uptake into the blood by meal ingestion, continuous enteral or intravenous nutrition. 
   \item Insulin production: $f_1$ represents the production of insulin. It is influenced by the concentration of glucose  with a delay $\tau_I$ to account for the time lag between high glucose levels triggering insulin production in the pancreas and when it becomes available for reducing glucose in the bloodstream. 
   \item Insulin-independent glucose utilization: $f_2$ describes the utilization of glucose by tissues, mainly the brain, in an insulin-independent manner. It does not rely on the presence of insulin.
   \item Insulin-dependent glucose utilization: $f_3\cdot f_4$ represents the utilization of glucose by muscle tissues in an insulin-dependent manner. It reflects the capacity of tissues to utilize insulin for glucose uptake. 
   \item Glucose production by the liver: $f_5$ represents the production of glucose by the liver. The delay $\tau_G$ represents the time between hepatic glucose production and insulin stimulation.
   \item Insulin degradation: The rate $d$ accounts for the degradation of insulin in the body, primarily by the liver and kidneys. It combines both natural factors (e.g., exercise) and artificial factors (e.g., medication) that influence the rate of insulin degradation. 
\end{compactitem}
The nonlinear pathways $f_1, f_2, f_3, f_4,$ and $f_5$ are represented using Hill functions, which are mathematical functions commonly used in biological modeling. 
These functions introduce additional parameters that have specific physiological interpretations and allow for a more accurate representation of the underlying dynamics of the glucose-insulin system, see methods section for details. 
%%%%%%%%%%%%%%%%%%%%%%%%%%%%%
The delays $\tau_I$ and $\tau_G$ are important physiological parameters encapsulating the responsitivity of the signaling and production pathways. 
They are assumed to be constant for the purpose of this article, although in practice, they can vary between individuals, as well as during the day and lifespan, and especially in the presence of diabetes.

% History of the model
The model has been extensively analyzed by various authors in the case of constant rates of glucose infusion \cite{li2006modeling,li2007analysis,huard2015investigation, huard2017mathematical}. 
It originates from the work of Sturis and collaborators who devised a model of glucose and insulin ultradian oscillations which were observed experimentally under various conditions \cite{sturisetal1991a}.
We also remark here that the model belongs to a larger class of models incorporating delays to capture secretion processes \cite{makroglou2006mathematical, shi2017oscillatory}.
We extend these earlier efforts on the analysis of the model by studying its response to periodic variations of the parameter $G_{\rm{in}}$, that is, periodic variations of glucose uptake. 
In particular, we consider on-off infusion, a form of periodic infusion that is comparatively easy to implement in practice, where the rate of glucose infusion periodically switches between a positive constant value and zero. 
Panels (b)--(e) of Fig.~\ref{fig:schema} show prototypical examples for the response of system \eqref{eq:G}--\eqref{eq:I} for various types of glucose uptake during fasting (b), glucose infusion with a (relatively high) constant rate (c) and periodic on-off infusion (d)--(e).  
We first investigate the loss of ultradian oscillations under sufficiently strong constant infusion, see Fig.~\ref{fig:schema}(b)--(c). 
We then aim to study the effects of different glucose infusion patterns $G_{\rm{in}}$ on glucose homeostasis, in particular the transiton from quasi-periodicity to entrainment, as shown in Figs.~\ref{fig:schema}(d)--(e). 

%in a scenario in which the system is subject to periodic perturbations in the form of repeated intravenous glucose intakes. %that temporarily increase the rate of food uptake. 
%The objective is to identify emerging dynamical structures which are shown to depend crucially on the period of the input. 
%is to study the level of variability and the phase of variability of the model output.

%
%While the physiology and modelling of intrapancreatic interactions have been successfully modelled and studied, the effects of inaccurate regulation on the overall feedback at the system level are not fully understood. 
%
%{\bf Circadian oscillation patterns in T1D\cite{vasquez2021oscillatory}}
% 
%In a similar manner, while the control of ultradian oscillations in glucose is typically impaired in individuals with decreased insulin sensitivity \cite{o1993lack}, oscillatory glucose infusions were shown to reintroduce a healthy rhythm \cite{sturis1993differential}.
%
 
%%Explain what is a typical test and how recent long-term test motivate investigating models with constant glucose infusion. 
%

%%%%%%%%%%%%%%%%%%
\section*{Results}
\subsection*{Ultradian oscillations}
It has been shown that for a fixed constant glucose infusion $G_{\rm{in}}$, sufficiently large values of the response delays $\tau_1$ and $\tau_2$ lead to periodic oscillations in system \eqref{eq:G}--\eqref{eq:I} with periods closely resembling the observed range for ultradian oscillations\cite{li2006modeling, huard2015investigation, huard2017mathematical}. 
Mathematically speaking, the onset of oscillations is mediated by a supercritical Hopf bifurcation that leads to a local topological change in the solution space of system \eqref{eq:G}--\eqref{eq:I} from a stable equilibrium to a situation of an unstable equilibrium surrounded by a small stable limit cycle close to the bifurcation point\cite{li2007analysis}. 
For details on bifurcation theory and the Hopf bifurcation, we refer the interested reader to Ref.~\citeonline{kuznetsov1998elements}. 
To witness the bifurcation point, it necessary to vary at least one parameter of the system. 
Here, we focus on the response delays $\tau_I$ and $\tau_G$. 
Allowing these two parameters values to vary simultaneously, one obtains a one-parameter curve $\mathbf{H(\omega)}=(\tau_I(\omega),\tau_G(\omega))$ of Hopf bifurcation in the $(\tau_I,\tau_G)$-plane in terms of the Hopf frequencies $\omega$, see Methods section for a detailed derivation. 
The curve $\mathbf{H(\omega)}$ corresponds to the critical curve for oscillations in system \eqref{eq:G}--\eqref{eq:I}. 

\begin{figure}[!]
    \centering
    \includegraphics[width=\linewidth]{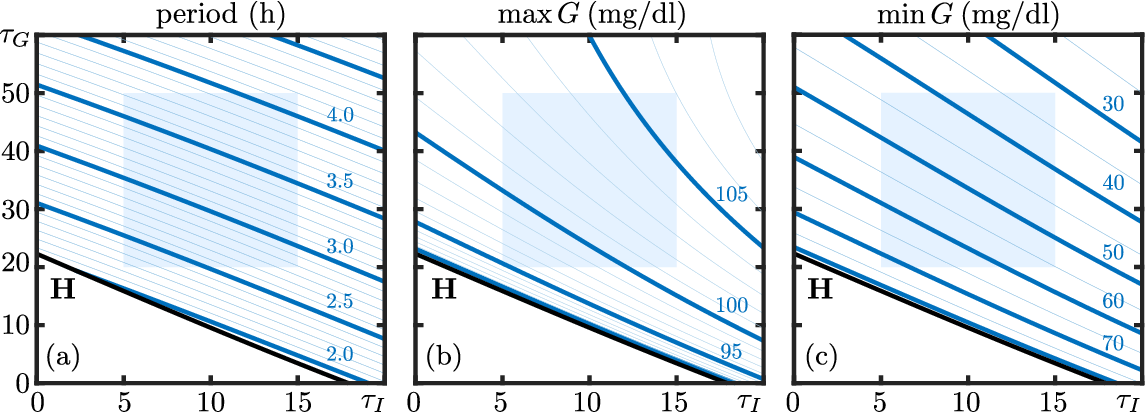}
    \caption{
    Characterization of fasting oscillations with respect to response delays. 
    Panels show the period (a), maximum glucose value (b), and minimum glucose value (c) as a function of response delays $\tau_I$ (min) and $\tau_G$ (min).
    Shown are the critical curve for oscillations (black, Hopf bifurcation), and iso-curves (blue) with constant period (a), glucose-maxima (b) and minimum of $G$ (c). 
    The light blue rectangle shows the physiological range of delay values for comparison. See methods section for the model and choice of parameters.}
    \label{fig:Hopfcurve-ext}
\end{figure}

Figure~\ref{fig:Hopfcurve-ext}(a) shows the curve $\mathbf{H}$ (black) during fasting, i.e. for $G_{\text{in}}=0$, computed with the software package DDE-Biftool for Matlab \cite{engelborghs2002numerical,sieber2014dde}.
It has been numerically verified that the curve $\mathbf{H}$ is indeed supercritical for the range of parameter values considered. 
Figure~\ref{fig:Hopfcurve-ext}(a) can be interpreted as follows: First, for value pairs above the curve and for $\tau_I\leq20~\mathrm{min},$ $\tau_G\leq60~\mathrm{min}$, any solution of the model starting in a physiological range of glucose and insulin develops periodic oscillations, see Fig.~\ref{fig:schema}(b).
Second, for value pairs $(\tau_I,\tau_G)$ below the curve $\mathbf{H}$, oscillations in system \eqref{eq:G}--\eqref{eq:I} decay and approach the equilibrium $(G^\ast,I^\ast).$
This possibly reflects the situation where an individual is administered a glucose dose that is too high to be managed in an oscillatory manner within physiological glucose and insulin ranges, compare Fig.~\ref{fig:schema}(c). 

Figure~\ref{fig:Hopfcurve-ext}(a) also gives an overview of the resulting period of oscillation above the critical curve $\mathbf{H}$ shown in the form of isocurves (blue) of limit cycles with constant period. 
The physiological range of parameters $(\tau_I,\tau_G)$ is highlighted by light blue square in the background for convenience. 
The range of expected periods for ultradian oscillations as predicted by the model thus ranges from $2.2$ to $4.2$ hours during fasting. 
More generally, we observe the period of the limit cycle oscillation grow approximately linear with the sum of the two delay values $\tau_I+\tau_G$. 
We also observe that away from the curve $\mathbf{H}$, the limit cycle oscillation becomes less and less sinusoidal, i.e. the nonlinearity of system \eqref{eq:G}--\eqref{eq:I} has more and more of an effect on the limit cycle. 
Panels (b)--(c) of Fig.~\ref{fig:Hopfcurve-ext} illustrate this effect by plotting isocurves of periodic orbit with constant minimum and maximum glucose within one period of oscillation. 
We observe that, whereas the glucose minimum decreases approximately linearly with the sum of the delays $\tau_I+\tau_G$, the maximum $G$ remains almost constant for the range of parameter values considered. 
Note that this predicted effect of long response delays is potentially harmful and is virtually undetectable by common testing methods.  

\begin{figure}[!]
    \includegraphics[width=0.7\linewidth]{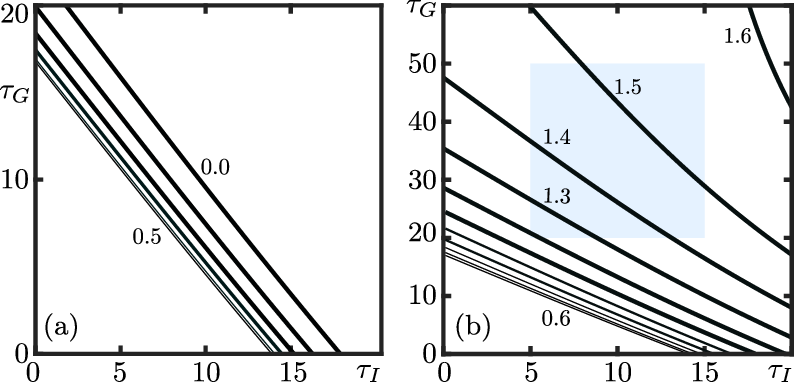}
    \caption{
    Position of the critical curve (curve of Hopf bifurcation) in the $(\tau_I,\tau_G)$-plane for various values of $G_{\text{in}}$ ranging from $0$ to $0.5$ and from $0.6$ to $1.6$ mg dl$^{-1}$ min$^{-1}$ (all black). 
    The light blue rectangle shows the physiological range of delay values for comparison.}
    \label{fig:Hopfcurve-GinVar}
\end{figure}

On the other hand, we observe that, for fixed values of the delays, gradually increasing the glucose infusion leads to a loss of oscillations. 
This phenomenon has been observed before and can be interpreted as an insufficient insulin secretion to accommodate the infusion, forcing the system to lower glucose lower levels\cite{li2007analysis}. 
Figure~\ref{fig:Hopfcurve-GinVar} shows how the location of the curve $\mathbf{H}$ changes for various levels of constant glucose infusion. 
We observe a two different types of change for values in the approximate ranges $0\leq G_\text{in}\leq 0.55$ mg dl$^{-1}$ min$^{-1}$ and $G_\text{in}>0.55$ mg dl$^{-1}$ min$^{-1}$ shown in panels (a) and (b) of Fig.~\ref{fig:Hopfcurve-GinVar}, respectively. 
Figure~\ref{fig:Hopfcurve-GinVar}(a) suggests that low levels of $G_\text{in}$ promote oscillations in system \eqref{eq:G}--\eqref{eq:I} as compared to the fasting case. 
This trend reverses at approximately at $G_\text{in}=0.55$ mg dl$^{-1}$ min$^{-1}$, where the location of the curve $\mathbf{H}$ starts moving to larger and larger values of $\tau_I$ and $\tau_G$, see Fig.~\ref{fig:Hopfcurve-GinVar}(b). 
Approximately at $G_\text{in}=1.2$ mg dl$^{-1}$ min$^{-1}$ the position of $\mathbf{H}$ is comparable with the starting location for $G_\text{in}=0$. 
Further increasing $G_\text{in}$ moves $\mathbf{H}$ inside the physiological range of delay values (light blue) and finally beyond causing all oscillations to cease in the physiological parameter regime. 
Compare also Fig.~\ref{fig:schema}(b)-(c) for an illustration of this transition and the loss of oscillations for $(\tau_I,\tau_G)=(5,20)$. 

\subsection*{Entrainment and amplitude response to on-off glucose infusion}
We now investigate the effect of periodic glucose infusion on baseline fasting oscillations shown in Fig.~\ref{fig:schema}(b), i.e. we fix $\tau_I=5$ min, $\tau_G=20$ min and periodically adjust the level of $G_\text{in}$ between $0$ and a positive value $G_{0}$ to be specified.
%This situation is analogous to performing repeated intravenous glucose tolerance tests, and it is interesting whether additional information can be obtained from such an infusion protocol. 
The natural frequency of ultradian oscillation in this case is $T\approx2.2$ h. 
We show that the resulting glucose and insulin ranges depend sensitively on the period of the on-off infusion.
Figures~\ref{fig:schema}(d)--(e) show two of the possible outcomes with different maximal infusion strength $G_{\max},$ period of infusion $T_\text{in}$, and infusion duration $t_{\text{in}}.$  

\begin{figure}[!]
    \centering
    \includegraphics[width=\linewidth]{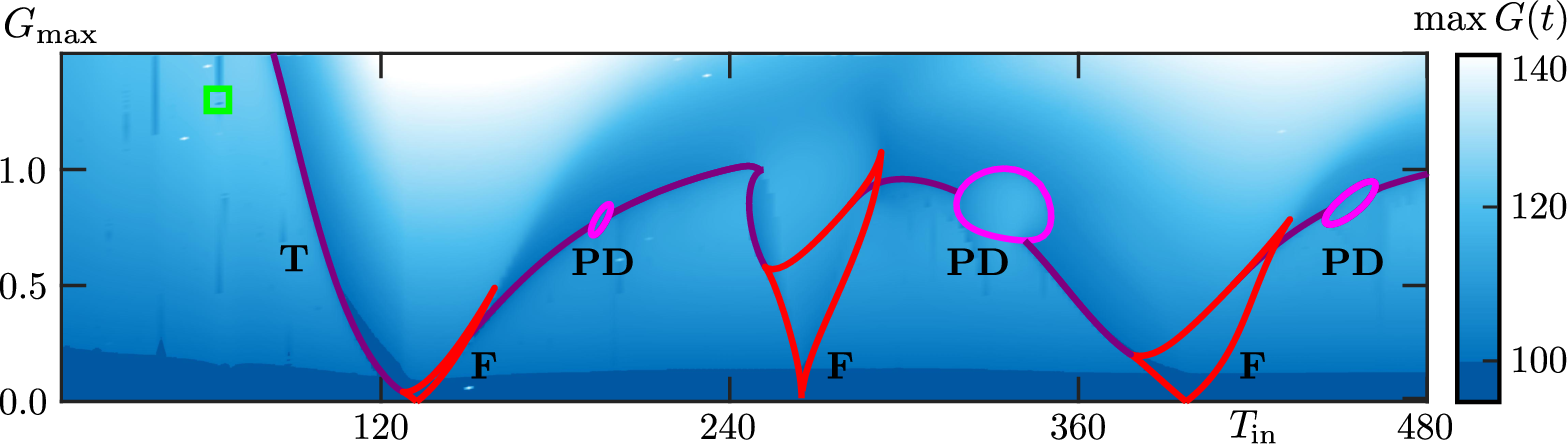}
    \caption{Response of model \eqref{eq:G}--\eqref{eq:I} to glucose infusion protocol \eqref{eq:Gin-periodic-forcing} with maximum infusion rate $G_{\max}$ (mg/(dl min)) and length of infusion $t_{\text{in}}=T_{\text{in}}/2$ (h). 
    Shown is the maximum value of $G$ (mg/dl) in colorcode (blue-white) obtained by integration for various $T_{\text{in}}$ and $G_{\max}$ over $100(\tau_I+\tau_G+T_{\text{in}})$ time units. 
    The maximum data is overlaid by curves of torus bifurcation (purple), curves of fold bifurcation of periodic orbits (red) and curves of period doubling bifurcation (magenta) bounding regions of locking to the infusion protocol. 
    Other parameters are $\tau_I=5$ min and $\tau_G=20$ min.}
    \label{fig:resonance}
\end{figure}

\paragraph*{Long infusion time compared to period}
Figure~\ref{fig:schema}(d) shows the result of periodic infusion with $G_\text{in}=1.35$ mg dl$^{-1}$ min$^{-1},$ for $t_\text{in}=30$ min every $T_\text{in}=60$ min, resulting in so-called quasi-periodic oscillations. 
Indeed, quasi-periodic oscillations are characterized by the presence of an oscillating envelope of the oscillation that evolves on a much slower time-scale, compare Fig.~\ref{fig:schema}(d). 
This is in sharp contrast with panels (b) (no infusion) and (c) (constant infusion with the same maximum rate) of Fig.~\ref{fig:schema}, where we have either periodic oscillations, or a decay of oscillations towards the equilibrium state. 
Quasi-periodic oscillations can be expected to occur in oscillatory systems which are externally driven by an input with non-commensurable period, here $T_\text{in}/T_0=2.2$. 
In this case, the effect of infusion very much depends on its current state: When insulin is low, glucose increases quickly; when insulin is high, glucose cannot increase further and the infusion only delays the expected decrease in glucose levels.

Periodicity of the oscillations can be restored by adjusting $G_\text{in}$ and  $T_\text{in}$. 
Figure \ref{fig:resonance} summarizes the response of system \eqref{eq:G}--\eqref{eq:I} to periodic forcing with different values of $G_\text{in}$ and  $T_\text{in}$. 
The locus in parameter space of the quasi-periodic oscillation shown in Fig.~\ref{fig:schema}(d) is indicated by a green rectangle. 
Figure \ref{fig:resonance} shows the overall glucose maximum (in color code) observed over a time span of $100\cdot(T_{\text{in}}+\tau_I+\tau_G)$ minutes. 
The various mechanisms generating periodic rhythms can be understood from the numerically computed bifurcation curves shown in \ref{fig:resonance}.
These correspond to curves of torus bifurcations $\mathbf{T}$ (purple), curves $\mathbf{F}$ (red) of fold (or saddle-node) bifurcations of periodic orbits, and curves $\mathbf{PD}$ (magenta) of period-doubling bifurcations of periodic orbits.
These mark the transition to periodic solutions and thus characterize the so-called \emph{entrainment} of oscillations.
%Additionally shown are numerically computed curves $\mathbf{T}$ (purple) of torus bifurcations, curves $\mathbf{F}$ (red) of fold (or saddle-node) bifurcations of periodic orbits, and curves $\mathbf{PD}$ (magenta) of period-doubling bifurcations of periodic orbits, indicating the transition to periodic solutions, so-called \emph{entrainment} of oscillations. 

The curves $\mathbf{F}$ respectively enclose deltoid-like regions -- called resonance or locking tongues -- extending from the line $G_{\text{in}}=0$, inside of which we observe periodic oscillations. 
The curves $\mathbf{F}$ emerge pairwise from resonant points where the infusion period is a rational multiple of the natural period of the system without infusion, i.e. $pT_{\text{in}}=qT_0$ for integers $p,q$. 
Figure~\ref{fig:resonance} shows the first three principal resonances of system \eqref{eq:G}--\eqref{eq:I} where $p=1,2,3$ and $q=1$. It is expected that such resonance tongues emanate from the line $G_{\max}$ at every rational point $T_0$. These higher order resonances (exept $p=4$ and $q=1$ which is outside of the considered range of parameter values) has been omitted/not computed as they are are typically very narrow and thus unlikely to be physiologically relevant.

This behavior persists moving towards larger values of $G_{\max}$ into the regions that are bounded approximately by the curves $T$, where the underlying stable periodic orbit destabilizes and gives rise to a torus that corresponds to quasi-periodic oscillations. 
We find numerical evidence that the direction with which this torus emanates from the the curve $\mathbf{T}$ can change and gives rise to the discontinuous transition between the observed maximum values in Fig.~\ref{fig:resonance}. 
The curves $\mathbf{T}$ each emanate from either point of intersection with a curve $\mathbf{F}$ or $\mathbf{PD}$. Intersections with curves $\mathbf{PD}$ correspond to higher order locking between the ultradian oscillations and the infusion. 
%A more detailed analysis of the stability of periodic orbits and description of the codimension-two bifurcation points is beyond the descriptive nature of the present article. 
%We simply mention that the curves $\mathbf{F},\mathbf{T},\mathbf{PD}$ mark the transition as they give rise to stable periodic oscillations. 
%The regions enclosed by these curves are commonly called resonance tongues or locking tongues. 
Overall, we observe that the strength and period of the infusion have a crucial effect on the resulting amplitude of the oscillations. 
For instance, forcing the system periodically with $T_0=T_\text{in}$ and relative amplitude $G_\text{in}=1$ mg dl$^{-1}$ min$^{-1}$ leads to a $40\%$ increase of the overall amplitude of the oscillation (which appears to be still in physiological range). 
In contrast, stimulating the system with a gradually increasing $G_\text{in}$ in the 2:1 regime first goes through phase during which glucose amplitudes remain relatively constant before slowly increasing. 

%amplitudes remain steady stimulating the system in the 2:1 regime leads 

More generally, we observe that, for the assumed values of the response delays, periodic infusion with $T_{\text{in}}=2t_{\text{in}}>T_0$ and $G_\text{in}$ is sufficient for the resulting period of the resulting glucose-insulin oscillation to be set by (locked to) the period of glucose infusion. 

%Up to three levels of \textbf{subheading} are permitted. Subheadings should not be numbered.

%We provide a numerical bifurcation analysis of the system with respect to this periodic forcing. 
%For a constant glucose infusion of $G_{in}=1.0 mg dl^{-1}$ and heuristic parameter values $\tau_1=5$ and $\tau_2=20$, we find an unstable glucose-insulin equilibrium and observe stable ultradian periodic oscillations with period $T=132.3$ (mins). 

%\begin{figure}[h]
%\begin{flushleft}
%\includegraphics[width=0.5\linewidth]{Fig/GI.png}\includegraphics[width=0.5\linewidth]{Fig/GI-limit-cycle.png}
%\end{flushleft} 
%\caption{Stable limit cycle in system \eqref{eq:G}-\eqref{eq:I} for a constant glucose infusion $G_{in}=1.0$ mg dl$^{-1}$ and delays $\tau_G=5$ and $\tau_I=20$.}
%\label{fig:phasespace}
%\end{figure}

%See Fig.~\ref{fig:phasespace} for a projection into the $(I,G)$-plane. 

\begin{figure}[!]
    \includegraphics[width=0.7\linewidth]{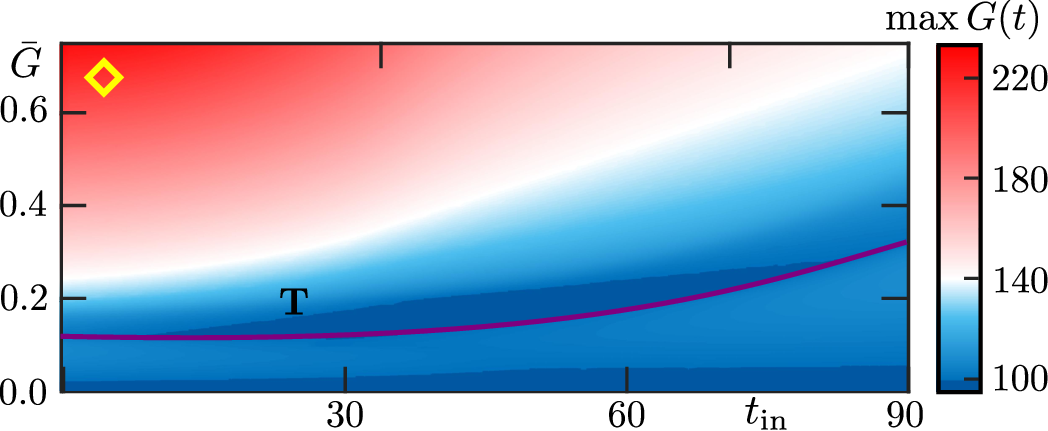}
    \caption{
    Response of model \eqref{eq:G}--\eqref{eq:I} to glucose infusion protocol \eqref{eq:Gin-periodic-forcing} with average infusion rate $\bar G=G_{\max} \cdot t_{\text{in}}/T_{\text{in}}$ mg dl$^{-1}$ min$^{-1}$ over the length of infusion $t_{\text{in}}$ (min) with constant period $T_{\text{in}}=180$ (min). 
    Shown is the maximum value of $G$ (mg/dl) in colorcode (blue-white) obtained by integration for various  $t_{\text{in}}$ and $ \bar G$. 
    The maximum data is overlaid by curves of torus bifurcation (purple), curves of fold bifurcation of periodic orbits (red) and curves of period doubling bifurcation (purple) bounding regions of locking to the infusion protocol. 
    Other parameters are $\tau_I=5$ min and $\tau_G=20$ min.}
    \label{fig:resonance2}
\end{figure}

\paragraph*{Short infusion time compared to period}
We note here that locking can be achieved when the same glucose dose is delivered in a shorter period of time, resulting in a more concentrated and intense infusion. To further explore this phenomenon, we conducted additional experiments using an on-off glucose infusion protocol with a fixed infusion period of $T_{\rm{in}}=180$ min. Figure~\ref{fig:resonance2} showcases the results obtained from these experiments, where we varied both the infusion time $t_{\rm{in}}$ and the average glucose dose per minute $\bar G$, represented by $G_{\max}\cdot t_{\rm{in}} / T_{\rm{in}}$.

In this figure, we observe a locus in the parameter space that corresponds to the quasi-periodic orbit illustrated in Fig.~\ref{fig:schema}(e). 
This locus is denoted by a distinctive yellow diamond marker, which highlights the specific combination of infusion time and glucose dose that leads to the observed quasi-periodic behavior. 
Additionally, we present a curve labeled as $\mathbf{T}$, which represents a torus bifurcation curve. 
This curve serves as an indicator of the critical transition point between entrainment and quasi-periodic oscillation in response to the infusion protocol.

%\paragraph*{Timing}
%We have explored the influence of the relative timing $\sigma_{\text{in}}$ between the natural oscillation and the glucose infusion, the co-called phase difference. A similar effect has been studied in Ref.~\ref{}. Contrary to Ref.~\ref{}, we have neither found a dependence of the overall glucose maximum, nor the position of the curves $\mathbf{F},\mathbf{T},\mathbf{PD}$ on the value of $\sigma_{\text{in}}$. Results not shown. 

\section*{Discussion}
%The Discussion should be succinct and must not contain subheadings.
It is well documented that glucose rhythms stimulate pulsatile pancreatic insulin secretion at various timescales \cite{sturisetal1991a,satin2015pulsatile}.
For example, the 1:1 entrainment mode -- namely one ultradian glucose oscillation per glucose infusion cycle -- was clinically shown to be present using a sinusoidal glucose infusion in individuals without diabetes \cite{o1993lack,sturis1995phase}.
Our analysis of periodically driven ultradian oscillations highlights that a periodic on-off stimulus, closer to normal daily conditions, also possesses the ability to entrain glucose rhythms.
Furthermore, the duration of each glucose input has a crucial impact on the generation of periodic rhythms, as well as on attained glycemic levels.
This theoretically provides a method for delivering a fixed glucose dose while minimising the amplitude of the resulting rhythm.
This can be achieved by either altering the period of the infusion, or the length of each pulse. 
This is most observable in figure \ref{fig:resonance2}, where stretching the infusion duration leads to lower glucose amplitudes. 
For example, consider a scenario where glucose is infused every 180 minute over a 12-hour period.
Infusing a dose with $G_{\max} = 2.4$ mg dl$^{-1}$ min$^{-1}$ over $t_{in}=30$ minutes leads to a maximal glucose value around $150$ mg/dl.
In contrast, a dose with $G_{\max} = 1.2$ mg dl$^{-1}$ min$^{-1}$ over $t_{in}=60$ minutes reduces the maximal glucose level to around $125$ mg/dl.
In both cases, the average dose per minute is $\bar{G} = 0.4$ mg dl$^{-1}$ min$^{-1}$, and a total dose of $288$ mg/dl is infused over the 12 hour timespan.

Our study provides valuable insights into the system's response to glucose infusion patterns, providing multiple pathways for the production of stable oscillatory rhythms and a similar entrainment structure is also expected for simpler models of glucose-insulin regulation featuring delays, e.g. \cite{panunzi_discrete_2007,li_range_2012,shi2017oscillatory}.
Nonetheless, there are several limitations that should be considered. 
First, let us note that while the exact location of bifurcation curves would depend on model parameters, the bifurcation types are likely to remain the same for parameter ranges representing non-diabetic individuals.
Our model assumes fixed values for the delays in insulin and glucose production pathways, represented by $\tau_I$ and $\tau_G$ respectively. 
%We also emphasise that while delays are assumed to be constant, only slight variations should be observed for smooth deviations of these values.
In reality, these delays can vary between individuals and change over short and long timescales due to daily-life factors such as exercise, aging and the presence of insulin resistance.
Future research could incorporate individual-specific delays to account for this variability and investigate their impact on the system's dynamics.

It is worth noting that our model relies solely on plasma glucose and insulin measurements for prediction, which highlights the importance of accurate and reliable measurements in clinical settings. 
The nonlinear structure of the model allows for the description of nontrivial dynamics and enhances parameter identifiability. 
This aspect is crucial for developing robust and accurate models that can capture the complex dynamics of the glucose-insulin regulatory system.
%Secretion times are likely to depend explicitly on glucose and insulin levels in a smooth manner.

It is also worth noting that the timing of the glucose infusion does not influence the bifurcation structure (\ref{fig:resonance}), nor the glucose-insulin ranges of the periodic rhythms.
In other words, the long term dynamics is not dependent on the starting time of the periodic on-off glucose infusion.
This does not mean that the timing of glucose inputs bears little importance. 
While the investigated infusion ranges ensured 
the positivity of glucose and insulin values, values below or above healthy physiological ranges may appear in the transient path to the limit cycle.
In turn, additional dynamics may emerge from interactions with other physiological feedback loops or subsystems, such as the glucagon pathway or the hypothalamic-pituitary-adrenal axis, for which the alignment with glucose regulation is essential for maintaining good health \cite{zavala2022misaligned}.
The recent incorporation of glucagon \cite{cohen2021novel} in models of the glucose-insulin feedback system may help provide a more complete and quantitative picture of dynamical interactions occurring within the pancreas \cite{Pedersen201329, montefusco2020heterogeneous} which can be used to improve quantitative tests for the detection and measurement of insulin and glucagon resistance \cite{morettini2021mathematical}. 

%It is important to acknowledge the limitations of our study. Our model assumes fixed values for the delays in insulin and glucose production pathways, which may vary between individuals and change over time. Incorporating individual-specific delays in future research could enhance the model's accuracy and capture inter-individual variability. Additionally, our model focuses solely on the glucose-insulin loop and does not consider the complex interactions between other physiological processes. Integrating these interactions into a more comprehensive model would provide a more realistic representation of the system's behavior.

Another aspect to consider is the interaction between the glucose-insulin regulatory system and other physiological processes. 
Our model focuses solely on the glucose-insulin loop, but in reality, there are complex interactions between various metabolic pathways, hormones, and organs. 
Integrating these interactions into a comprehensive model could provide a more complete understanding of the system's behavior and its response to different stimuli.

\section*{Conclusion}
In this study, we employed a system-level mathematical model to investigate the response of the glucose-insulin regulatory system to periodic glucose infusion. 
By exploring different glucose infusion patterns and analyzing the resulting dynamics, we gained insights into the system's behavior and identified key factors influencing its response.

Our findings demonstrate that the glucose-insulin regulatory system exhibits a range of behaviors depending on the glucose infusion pattern. When a constant glucose infusion is applied, the system shows ultradian oscillations characterized by periodic variations in glucose and insulin levels. However, as the glucose infusion rate exceeds a certain threshold, these oscillations disappear, and the system focuses on reducing glucose levels without exhibiting oscillatory behavior. This observation suggests a physiological limit beyond which the system's oscillatory capacity is overwhelmed.

We further investigated the effects of periodic on-off pulses, mimicking repeated intravenous glucose tolerance tests. Our analysis revealed that the period of the on-off pulses plays a crucial role in determining the system's dynamics and glucose-insulin ranges. Different patterns of oscillations, including stable limit cycles and irregular oscillations, were observed for varying infusion periods. This highlights the importance of considering the frequency and duration of glucose stimuli in understanding the system's response.

The results of this study have important implications for understanding glucose regulation in both normal and abnormal physiological conditions. 
By elucidating the system's response to different glucose infusion patterns, our findings can inform the development of test strategies for evaluating the system's performance and identifying potential dysfunctions. 
Furthermore, they provide insights into the underlying mechanisms governing glucose-insulin dynamics, contributing to the broader understanding of metabolic regulation.

In conclusion, our study enhances our understanding of the glucose-insulin regulatory system by investigating its response to periodic glucose infusion. 
We identified the impact of different glucose infusion patterns on the system's dynamics and demonstrated the importance of various types of glucose stimuli. 
These insights can aid in the development of diagnostic and therapeutic strategies for glucose regulation and contribute to advancements in the management of metabolic disorders. 
Future research should aim to incorporate individual-specific delays and consider the broader physiological context to further refine our understanding of glucose regulation and its implications for human health.

\section*{Methods}

\subsection*{The glucose-insulin regulatory delayed-feedback model}

We consider the system-level mathematical model 
\begin{align}
    G'(t) &= G_{\text{in}}(t) - f_2(G(t)) - f_3(G(t))f_4(I(t)) + f_5(I(t-\tau_G)) \label{eq:G}\\
    I'(t) &= I_{\text{in}}(t) + f_1(G(t-\tau_I)) - d I(t) \label{eq:I}
\end{align}
with variables $I(t)$ and $G(t)$ representing the concentrations  and  of glucose (mg dl$^{-1}$) and insulin (uU ml$^{-1}$) in the plasma at time instant $t$. 
System (\ref{eq:G})--(\ref{eq:I}) explicitly depends on time delays $\tau_I$ and $\tau_G$ respectively representing the system's response time to insulin production as a result to glucose uptake, and the production of glucose by the liver as a result of low insulin levels. 
Glucose intake and insulin infusion are modeled by parameters, here called $G_{\rm{in}}$ and $I_{\rm{in}}$. The physiological response of body is modeled by the nonlinearities
\begin{align*}
    f_1(G) &= \frac{R_m G^{h_1}}{G^{h_1} + (V_g k_1)^{h_1}},\\
    f_2(G) &= \frac{U_b G^{h_2}}{G^{h_2} + (V_g k_2)^{h_2}},\\
    f_3(G) &= \frac{G}{C_3 V_g},\\
    f_4(I) &= U_0 + \frac{ (U_m - U_0)I^{h_4}}{I^{h_4} + (1/V_i + 1/(Et_i))^{-h_4}k_4^{h_4}},\\
    f_5(I) &= \frac{R_g I^{h_5}}{I^{h_5} + (V_p k_5)^{h_5}},
\end{align*}
where $R_m=210,$ $V_i=11$, $V_g=10$, $E=0.2$, $U_b=72$, $t_i=100$, $C_3=1000$, $R_g=180$, $U_0=40$,
$V_p=3$, $U_m=940$, $h_1=2$, $k_1=6000$, $h_2 = 1.8$, $k_2=103.5$, $h_4 = 1.5$, $k_4=80$, $h_5=-8.54$, and $k_5=26.7$ with corresponding units. Insulin degradation is modeled by a constant rate $d$. Throughout the paper we fix $d=0.06$. 
The model has been considered before and has been analyzed extensively by various authors \cite{li2006modeling,li2007analysis,huard2015investigation, huard2017mathematical}. In particular, it can be shown that, for the parameter values considered and in the absence of infusion, there is a unique equilibrium solution $(G^\ast,I^\ast)$; see for example [\citeonline{bennett2004global}].
The delay parameters used for numerical simulation are $\tau_I=5$ and $\tau_G=20$ if not stated otherwise. 
For the general theory of delay differential equations, such as existence, uniqueness and the stability of solutions, we refer the interested reader to classic textbooks on the topic\cite{hale2013introduction,diekmann2012delay}. 

\subsection*{Critical delay values for oscillatory behavior when infusion rate is constant}
The critical curve for oscillations in the ($\tau_I,\tau_G$)-parameter plane can be computed from the linearization of system (\ref{eq:G})--(\ref{eq:I}) about the equilibrium solution $(G^\ast,I^\ast)$ and imposing the condition $\lambda=i\omega,~\omega>0$ (Hopf bifurcation) on solutions of the corresponding characteristic equation
\begin{equation}\label{eq:char}
    0=\chi(\lambda):=\lambda^2+\alpha_1 \lambda +\alpha_0
+\beta_1 e^{-\lambda \tau_1} + \beta_2 e^{-\lambda \tau_2}, 
\end{equation}
where $\tau_1=\tau_I,$ $\tau_2=\tau_I+\tau_G$ and  $\alpha_1=f_2^\prime
(G^\ast) + f^\prime_3(G^\ast)f_4(I^\ast)+d$, $\alpha_0=d(f_2^\prime
(G^\ast) + f^\prime_3(G^\ast)f_4(I^\ast))$, $\beta_1=f^\prime_
1(G^\ast)f_3(G^\ast)f^\prime_4
(I^\ast)$, $\beta_2=-f^\prime_1
(G^\ast)f^\prime_5(I^\ast)$.  A detailed derivation of Eq.~(\ref{eq:char}) can be found in [\citeonline{huard2017mathematical}].\\
The equation $0=\chi(i\omega)$ can be solved parametrically for $\tau_1$ and $\tau_2$ to give
\begin{align}
\tau_{1,2}(\omega)&=\frac{1}{\omega}\left(\arctan\left(\frac{\alpha_1\omega}{\omega^2-\alpha_0}\right) + \arccos\left(\frac{\beta_{2,1}^2-\beta_{1,2}^2-(\omega^2-\alpha_0)^2 -\alpha_1^2\omega^2}{2\beta_{1,2}\sqrt{(\omega^2-\alpha_0)^2 +\alpha_1^2\omega^2}}\right)
\right),\label{eq:tau_om}
\end{align}
revealing the critical curve for oscillations $\mathbf{H}\subset \mathbb{R}^2$ (curve of Hopf bifurcation)
\begin{align}
\mathbf{H}(\omega)&=(\tau_I(\omega),\tau_G(\omega))=(\tau_1(\omega),\tau_2(\omega)-\tau_1(\omega)).
\end{align} 
For the considered parameter values, we have that $\alpha_1 > \alpha_0$ and $\beta_2>\alpha_0$, ensuring the existence of $\mathbf{H}$. Indeed the curve is a sharp threshold for oscillation, as it can been shown numerically that for positive values ($\tau_I,\tau_G$) below $\mathbf{H}$ the fixed point $(G^\ast,I^\ast)$ is stable for any physiological range of starting values $G$ and $I$. 
It is worth noting here that system (\ref{eq:G})--(\ref{eq:I}) undergoes further Hopf bifurcations, respectively at $\tau_{I,k}(\omega)=\tau_{I}(\omega)+2\pi k/\omega$ and $\tau_{G,l}(\omega)=\tau_{G}(\omega)+2\pi l/\omega$ with $k,l$ an integer; however, for the parameter values considered, we can restrict ourselves to the smallest positive such value pair to cover the physiological parameter range. The range of relevant values of $\omega$ resulting in positive delays cannot be computed explicitly, however, straightforward calculations show that the boundaries $\omega_{I},\omega_{G}$ satisfying $\tau_I(\omega_{I})=0$ and $\tau_G(\omega_{G})=0$ are given by 
\begin{align}
\omega_G &= \sqrt{\alpha_0-\frac{\alpha_1^2}{2}
+\sqrt{\left(\alpha_0-\frac{\alpha_1^2}{2}\right)^2+(\beta_1+\beta_2)^2-\alpha_0^2}},\\
\omega_I &= \sqrt{\alpha_0+\beta_1-\frac{\alpha_1^2}{2}
+\sqrt{\left(\alpha_0+\beta_1-\frac{\alpha_1^2}{2}\right)^2+\beta_2^2-\alpha_0^2}},
\end{align}
with the corresponding delay values
\begin{align}
\tau_I(\omega_{G}) &=\frac{1}{\omega_{G}}\arctan \left(\frac{\alpha_1\omega_{G}}{\omega_{G}^2-\alpha_0}\right) + \frac{2\pi k^\ast}{\omega_G},\\
\tau_G(\omega_{I}) &=\frac{1}{\omega_{I}}\arctan \left(\frac{\alpha_1\omega_{I}}{\omega_{I}^2-\alpha_0-\beta_1}\right)+ \frac{2\pi l^\ast}{\omega_I},
\end{align}
where $k^\ast,l^\ast$ are the smallest integers such that $\tau_G$ and $\tau_I$ are positive.

We remark that the curve $\mathbf{H}$ vaguely resembles a straight line with slope $-1$ in the $(\tau_I,\tau_G)$-plane. This can be understood by exploiting the fact that parameter $\beta_1$ is small on the order of $10^{-3}$. Imposing the regular perturbation ansatz $\omega=\omega_0+\beta_1\omega_1+\mathcal{O}(\beta_1^2)$ on the imaginary part of (\ref{eq:char}) and comparing at zeroth and first order in $\beta_1$, we formally obtain
\begin{align}
\omega_0 &= \sqrt{\alpha_0-\frac{\alpha_1^2}{2}
+\sqrt{\left(\alpha_0-\frac{\alpha_1^2}{2}\right)^2+\beta_2^2-\alpha_0^2}},\\
\omega_1 &= \frac{\omega_0}{\alpha_1-\alpha_0+\omega_0^2}\tau_I \leq \frac{1}{2 \sqrt{\alpha_1-\alpha_0}} \tau_I.%
\end{align}

Thus, we can approximate $\mathbf{H}$ to first order in $\beta_1$ 
\begin{align}
\mathbf{H}(\omega_0 + \beta_1\omega_1(\tau_I))&\approx\left(\tau_I,\tau_2(\omega_0 + \beta_1\omega_1(\tau_I))-\tau_I\right)
\end{align}
by using the expression $\tau_2(\omega)$ in (\ref{eq:tau_om}).
As a result, $\mathbf{H}$ approaches the graph of the function $\tau_I\mapsto \tau_2(\omega_0)-\tau_I$ with slope $-1$ as $\beta_1\to 0$, which can be considered as a zero order approximation of $\mathbf{H}$.

\subsection*{On-off periodic infusion}
\begin{figure}[!]
\centering
    \includegraphics[width=0.75\linewidth]{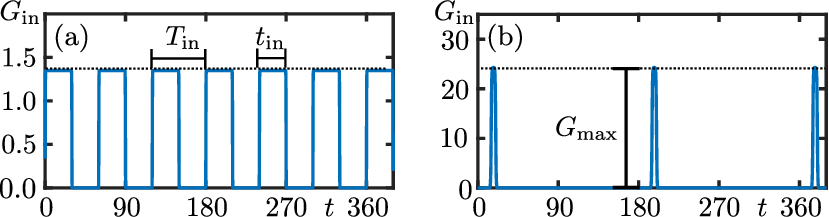}
    \caption{Form of glucose infusion used to obtain Figs.~\ref{fig:schema}(d)--(e).}
    \label{fig:forcing}
\end{figure}

Periodic shot infusion is modeled by a smooth, periodic, quickly varying function between $G_{\text{in}}(t)=I_{\text{in}}(t)=0$ in mg/(dl$\cdot$min) (no infusion) and $G_{\text{in}}(t)=G_{\max}$ (glucose infusion), 
%$I_{\text{in}}(t)=I_{\max}$ (insulin infusion)
respectively. We consider the specific form
\begin{equation}
\begin{split}
 G_{\text{in}}(t) = G_{\max}\cdot s(t-\sigma_G), \,
 %\quad  I_{\text{in}}(t) = I_{\max}\cdot s(t-\sigma_I),\\
 s(t)=h(\sin(2\pi t/T_{\text{in}}))\cdot h(\sin(2\pi(t-t_{\text{in}})/ T_{\text{in}} - \pi)),
 \end{split}
 \label{eq:Gin-periodic-forcing}
 \end{equation}
where the sigmoidal function $h(y)=(1+\exp(-k y))^{-1}$ can be considered a smooth version of the Heaviside step function $H(y)=0$ if $y<0$, and $H(y)=1$ if $y\geq0$ for sufficiently large $k$. The form of (\ref{eq:Gin-periodic-forcing}) was inspired by a model for auditory perception\cite{ferrario2021auditory}.
$T_{\text{in}}$ is the time between consecutive shots of glucose/insulin with duration $t_{\text{in}}$, and the lag $\sigma_G$ can be used to specify the timing of the infusion with respect to the underlying oscillation. 
The parameter $k$ models the initial and terminal variations at the beginning and end of the shot application. 
A detailed analysis of the influence of the $k$-parameter is beyond the scope of this paper. 
We found the choice $k=100$ to be sufficient. Panels (a) and (b) of Figure \ref{fig:forcing}  show the shape of the infusion patterns used for numerical computation of time series in Fig.~\ref{fig:schema}.

\subsection*{Numerical bifurcation analysis of time-delay systems with periodic infusion}

Numerical simulations have been obtained using pydelay \cite{flunkert2009pydelay}. 
Numerical bifurcation analysis has been performed using the software package DDE-BIFTOOL for Matlab/Octave\cite{sieber2014dde}. 
For a general introduction to numerical continuation methods available for delay differential equations and their application to physiological systems see Refs.~\citeonline{krauskopf2022bifurcation} and \cite{engelborghs2002numerical}, respectively. 
Isocurves in Fig.~\ref{fig:Hopfcurve-ext} have been computed using numerical continuation of periodic orbits in two parameters with the additional condition fixed period (a), and fixed maximum value (b), fixed maximum value (c), where in cases (b) and (c) we also relaxed the phase condition. 
For bifurcation analysis in the presence of periodic infusion, we append the two-dimensional ordinary differential equation
\begin{align}
    x^\prime (t) &= x - \omega y(t) - x(t)(x(t)^2+y(t)^2),\\
    y^\prime (t) &= -\omega x(t) + y  - y(t)(x(t)^2+y(t)^2),
\end{align}
with known stable periodic solution $(x(t),y(t))=(\cos(\omega t),\sin(\omega t))$ to system~(\ref{eq:G})--(\ref{eq:I}). The method has been employed in several other works, see for example Ref.~\citeonline{keane2015delayed}.
We achieve the specific form of infusion (\ref{eq:Gin-periodic-forcing}) by setting $G_{\text{in}}(t) = G_{\max} h(y(t-\sigma_G))h(y(t-\sigma_G-t_{in}\omega-\pi))$,
%and  $I_{\text{in}}(t) = I_{\max} h(y(t-\sigma_I))h(y(t-\sigma_I-\sigma_{\text{in}}))$
where $\omega=2\pi/T_{\text{in}}$.% and $ \sigma_{\text{in}} = \pi(1+2t_{\text{in}}/T_{\text{in}})$. 

\bibliography{main}
%\noindent LaTeX formats citations and references automatically using the bibliography records in your .bib file, which you can edit via the project menu. Use the cite command for an inline citation, e.g.  \cite{Hao:gidmaps:2014}.
%For data citations of datasets uploaded to e.g. \emph{figshare}, please use the \verb|howpublished| option in the bib entry to specify the platform and the link, as in the \verb|Hao:gidmaps:2014| example in the sample bibliography file.

\section*{Acknowledgements (not compulsory)}

The authors thank Jan Sieber for helpful discussions on the implementation of the numerical continuation methods in DDE-BIFTOOL.

\section*{Author contributions statement}
S.R. and B.H contributed equally to the study. Both authors conducted numerical experiments, analyzed the results, and contributed to the draft of the manuscript. Both authors reviewed the manuscript. 

\section*{Data available statement}
Numerical procedures to generate figures are available from the corresponding author on reasonable request. 

\section*{Additional information}
The authors declare no competing interests.

\end{document}